\documentclass[12pt]{article}
\usepackage{float}
\usepackage{amssymb}
\usepackage{amsthm}
\usepackage{amsmath}
\usepackage{stmaryrd}
\usepackage{titletoc}
\usepackage{mathrsfs}
\usepackage{authblk}
\usepackage[colorlinks, allcolors=blue]{hyperref}

\usepackage{amsfonts}
\usepackage{geometry}
\usepackage{graphicx}
\usepackage{float}
\usepackage{cite}
\usepackage[linesnumbered, ruled, vlined]{algorithm2e}

\newlength{\defbaselineskip}
\newcommand{\setlinespacing}[1]%
           {\setlength{\baselineskip}{#1 \defbaselineskip}}

\theoremstyle{plain}
\newtheorem{thm}{Theorem}[section]

\newtheorem{lem}[thm]{Lemma}

\theoremstyle{definition}

\newtheorem{ass}{Assumption}
\newtheorem{rem}{Remark}[section]
\newtheorem{ex}{Example}
\newcommand{\bR}{\mathbb{R}}

\newcommand{\bZ}{\mathbb{Z}}
\newcommand{\cA}{\mathcal{A}}

\newcommand{\cL}{\mathcal{L}}

\newcommand*{\dif}{\mathop{}\!\mathrm{d}}
\newcommand*{\E}{\mathbb{E}}
\newcommand*{\R}{\mathbb{R}}
\newcommand*{\J}{\mathcal{J}}

\makeatletter\@addtoreset{equation}{section} \makeatother

\title{Convergence of the Deep BSDE method for FBSDEs with non-Lipschitz coefficients}
\author[a]{Yifan Jiang \thanks{Email: yifan.jiang@maths.ox.ac.uk}}
\affil[a]{Mathematical Institute, University of Oxford}
\author[b]{Jinfeng Li \thanks{ Email: lijinfeng@fudan.edu.cn}}
\affil[b]{School of Mathematical Sciences, Fudan University}

\date{}

\begin{document}

\maketitle


\begin{abstract}
    This paper is dedicated to solving high-dimensional coupled FBSDEs with non-Lipschitz diffusion coefficients numerically.
    Under mild conditions, we provided a posterior estimate of the numerical solution that holds for any time duration.
    This posterior estimate validates the convergence of the recently proposed Deep BSDE method.
  In addition, we developed a numerical scheme based on the Deep BSDE method and presented numerical examples in financial markets to demonstrate the high performance.
  \end{abstract}


\section{Introduction}
Let \((\Omega,\mathcal{F},P)\) be a complete probability space, \(W\) a standard \(d\)-dimensional Brownian motion, \(T>0\) a fixed terminal time, \(\mathbf{F}:=\left\{ \mathcal{F}_{t} \right\}_{0\leq t\leq T} \) the natural filtration generated by \(W\) and augmented by the \(P\)-null sets.
By \(\cL_{\mathbf{F}}^{2}([0,T],\bR^{n})\) we denote the space of \(\mathbf{F}\)-adpated \(\bR^{n}\)-valued stochastic process \((X_{t})_{0\leq t\leq T}\) satisfying \(\E\int_{0}^{T}|X_{t}|^{2}\dif t<\infty\).
 The following coupled forward-backward stochastic differential equation (FBSDE) is considered in this paper:
\begin{equation}
    \label{eqn-fbsde}
    \left\{
    \begin{aligned}
        X_{t} & = x + \int_{0}^{t}b(s,X_{s},Y_{s})\dif s+ \int_{0}^{t}\langle\sigma(s,X_{s}),\dif W_{s}\rangle, \\
        Y_{t} & =g(X_{T})+\int_{t}^{T}f(s,X_{s},Y_{s},Z_{s})\dif s-\int_{t}^{T}\langle Z_{s},\dif W_{s}\rangle,
    \end{aligned}
    \right.
\end{equation}
where \(\left( X_{t} \right)_{0\leq t\leq T}\), \(\left( Y_{t} \right)_{0\leq t\leq T}\), \(\left(Z_{t} \right)_{0\leq t\leq T}\) are in \(\cL_{\mathbf{F}}^{2}([0,T],\bR^{n})\), \(\cL_{\mathbf{F}}^{2}([0,T],\bR^{m})\), \(\cL_{\mathbf{F}}^{2}([0,T],\bR^{m\times d})\) respectively.

The coefficients \(b\), \(\sigma\), \(f\), \(g\) are all deterministic.
The existence and uniqueness of the solution to the fully-coupled FBSDEs has been widely explored and we refer to \cite{peng1991probabilistic,ma1994solving,peng1999fully,pardoux1999forward,ma1999forward} for details.
Recently, the Deep BSDE method \cite{weinan2017deep} was proposed for numerically solving high-dimensional BSDEs and parabolic PDEs.
This method has high efficiency in approximating the accurate solution due to neural networks' universal approximation capability.
The convergence of the Deep BSDE method for FBSDEs has been extensively studied when the coefficients are sufficiently regular, see \cite{han2020convergence, ji2020three}.
A posterior estimate was used to bound the estimated error of the numerical solution.
Hence, the original numerical solution problem can be reformulated as a stochastic optimization problem.
A neural network is then used to optimize the posterior estimate and obtain the numerical solution.

In this paper, we apply the Deep BSDE method to FBSDEs with non-Lipschitz diffusion coefficients which have important applications, e.g., Cox-Ingersoll-Ross (CIR) model \cite{cox2005theory}.
The forward stochastic differential equation admits a unique solution when the diffusion coefficient \(\sigma\) is uniformly H\"older continuous with order \(\frac{1}{2}\).
The detailed proof can be found in \cite{yamada1971uniqueness,yamada1971uniqueness2}.
The CIR process has received a lot of attention in computational finance.
The strong convergence of various discretization schemes has been demonstrated in \cite{deelstra1998convergence, DereichCIR, Alfonsi2015}.
Our main contribution is to provide a posterior estimate for non-Lipschitz FBSDEs that is valid for any time duration.
To the best of our knowledge, this is the first theoretical result that supports the convergence of the Deep BSDE method for non-Lipschitz FBSDEs.
Due to the non-Lipschitz diffusion coefficient, it is difficult to balance the order of the estimate between the forward and backward equations.
We apply a series of Yamada-Watanabe functions \cite{yamada1971uniqueness,yamada1971uniqueness2} to solve this imbalance and bound the error of the time discretization.
In the previous studies \cite{han2020convergence,ji2020three}, the diffusion coefficient \(\sigma\) of the forward equation is assumed to be Lipschitz continuous and the time duration is required to be sufficiently small.
We extend the posterior estimate to arbitrary time duration only under mild assumptions on the decoupling field of FBSDE \eqref{eqn-fbsde}.

This paper will be organized as follows.
In the following section, we first formulate our main results, including theorems and numerical algorithms. Sections 3–5 are devoted to proving the main theorems.
Section 6 includes numerical examples to demonstrate the application in finance.
\section{Main Results}
\subsection{Assumptions}
The problem we mainly focus on is the numerical algorithm for FBSDEs with non-Lipschitz diffusion coefficient \(\sigma\).
We will make the following assumptions:
\begin{ass}
    We assume the diffusion coefficient \(\sigma:[0,T]\times\R^{n}\to\R^{n\times d} \) has the following form
    \[\sigma(t,x)^{\intercal}=\left( \sigma_{1}(t,x_{1})^{\intercal}, \cdots, \sigma_{n}(t,x_{n})^{\intercal}\right),\]
    where \(x_{i}\) is the \(i\)-th component of \(x\).
\end{ass}

\begin{ass}
    Coefficients \(b,\, f,\,g\) are uniformly Lipschitz continuous and \(\sigma_{i}\) is uniformly H\"older continuous with order \(\frac{1}{2}\).
    \begin{enumerate}
        \item \(|b(t,x,y)-b(t',x',y')|\leq L_{b}(|t-t'|+|x-x'|+|y-y'|),\)
        \item \(|f(t,x,y,z)-f(t',x',y',z')|\leq L_{f}(|t-t'|+|x-x'|+|y-y'|+|z-z'|),\)
        \item \(|g(x)-g(x')|\leq L_{g}|x-x'|,\)
        \item \(|\sigma_{i}(t,x_{i})-\sigma_{i}(t',x_{i}')|^{2}\leq L_{\sigma}(|t-t'|+|x_{i}-x_{i}'|),\)
    \end{enumerate}
    where \(L_{b}\), \(L_{f}\), \(L_{g}\), and \(L_{\sigma}\) are given positive constants.
\end{ass}

\begin{rem}
    We emphasize that the specific structure of \(\sigma\) is required to use the Yamada-Watanabe trick in \cite{yamada1971uniqueness}.
    The uniqueness of the solution may fail for general \(\sigma\) according to \cite[Remark 2]{yamada1971uniqueness2}.
\end{rem}

It is well-known that FBSDE \eqref{eqn-fbsde} is closely related to a parabolic PDE \cite{peng1991probabilistic, ma1994solving}.
\begin{ass}
    We assume the following PDE 
    \begin{equation}
        \label{eqn-decouple}
        \left\{
        \begin{aligned}
             & u_{t}+\nabla_{x} ub(t,x,u)+\frac{1}{2}\langle  \sigma^{\intercal}\sigma(t,x), \mathbf{H}_{x} u\rangle+f(t,x,u,\nabla_{x} u\sigma(t,x))=0, \\
             & u(T,x)=g(x).
        \end{aligned}
        \right.
    \end{equation}
    has a unique classical solution \(u(t,x):[0,T]\times\R^{n}\to\R^{m}\).
    Here \(\mathbf{H}_{x}\) denotes the Hessian with respect to argument \(x\).
\end{ass}

\begin{ass}
    We assume \(u\) is uniformly Lipschitz continuous,
    \begin{equation*}
        |u(t,x)-u(t,x')|\leq L_{u}|x-x'|,
    \end{equation*}
    where \(L_{u}\) is a given positive constant.
\end{ass}
\begin{rem}
    Equation \eqref{eqn-fbsde} has a unique solution if Assumptions 1-4 hold, see \cite{ma1994solving,zhang2017backward}.
    \(u\) is called the decoupling field of FBSDE \eqref{eqn-fbsde}.
    In particular, the solution of equation \eqref{eqn-fbsde} satisfies \(Y_{t}=u(t,X_{t})\) and \(Z_{t}=\nabla_{x}u(t,X_{t})\sigma(t,X_{t})\).
\end{rem}

\subsection{Main Theorems}
To solve FBSDEs \eqref{eqn-fbsde} numerically, we introduce the following stochastic optimal control problem.
Consider an \((n+m)\)-dimensional SDE system \((\tilde{X},\tilde{Y})\):
\begin{equation}
    \label{eqn-fsde}
    \left\{
    \begin{aligned}
        \tilde{X}_{t} & = x + \int_{0}^{t}b(s,\tilde{X}_{s},\tilde{Y}_{s})\dif s+ \int_{0}^{t}\langle\sigma(s,\tilde{X}_{s}),\dif W_{s}\rangle,         \\
        \tilde{Y}_{t} & =\tilde{y}-\int_{0}^{t}f(s,\tilde{X}_{s},\tilde{Y}_{s},\tilde{Z}_{s})\dif s+\int_{0}^{t}\langle\tilde{Z}_{s},\dif W_{s}\rangle,
    \end{aligned}
    \right.
\end{equation}
with the aim of minimizing the objective functional
\begin{equation}
     \J(\tilde{y},\tilde{Z}_{.})=\left(\E|g(\tilde{X}_{T})-\tilde{Y}_{T}|^{2}\right)^{\frac{1}{2}}.
\end{equation}
We view \((\tilde{y},\tilde{Z}_{.})\) as the control of the system and assume \(\tilde{Z}_{.}\in \cL_{\mathbf{F}}^{2}([0,T],\bR^{m\times d})\). 
Under the regularity assumption on the coefficients, we have the existence and uniqueness of system \eqref{eqn-fsde}.
This result can be found in \cite{yamada1971uniqueness,yamada1971uniqueness2}.
It is simple to verify that \((y,Z_{.})\) the solution of FBSDE \eqref{eqn-fbsde} is the optimal control.

Our first result demonstrates that the difference between any control and the optimal control can be bounded by the objective functional.
Thus, we can solve the original FBSDE \eqref{eqn-fbsde} by finding the optimal control of system \eqref{eqn-fsde}.

\begin{thm}
    \label{thm-control}
    We assume Assumptions 1-4 hold.  
    Let \((X_{t},Y_{t})\) be the solution to FBSDE \eqref{eqn-fbsde} and \((\tilde{X}_{t},\tilde{Y}_{t})\) be the state of system \eqref{eqn-fsde} under the control \((\tilde{y},\tilde{Z}_{.})\in \bR^{m}\times \cL^{2}_{\mathbf{F}}([0,T],\bR^{m\times d})\).
    We have for some constant \(C\) 
    \[\sup_{0\leq t\leq T}\E|X_{t}-\tilde{X}_{t}|+\sup_{0\leq t\leq T}\E|Y_{t}-\tilde{Y}_{t}|\leq C  \J(\tilde{y},\tilde{Z}_{.}).\]
    Furthermore, we have 
        \[\sup_{0\leq t\leq T}\E|X_{t}-\tilde{X}_{t}|^{2}+\sup_{0\leq t\leq T}\E|Y_{t}-\tilde{Y}_{t}|^{2}+\int_{0}^{T}\E|Z_{t}-\tilde{Z}_{t}|^{2}\dif t\leq C\left[  \J(\tilde{y},\tilde{Z}_{.})+ \J(\tilde{y},\tilde{Z}_{.})^{2}\right].\]
    Here the constant \(C\) does not depend on the choice of \((\tilde{y},\tilde{Z}_{.})\).
\end{thm}

Then, we consider the Euler-Maruyama scheme of system \eqref{eqn-fsde} to obtain the numerical solution to equation \eqref{eqn-fbsde}.
Let \(\Delta t=TN^{-1}\) be the discrete-time step and \(\tau(t)=\lfloor t\Delta t^{-1}\rfloor \Delta t\).
\((\hat{X}_{t},\hat{Y}_{t})\) is the state of the following discrete system under the control \((\tilde{y},\tilde{Z}_{.})\):
\begin{equation}
    \label{eqn-dis-fsde}
    \left\{
    \begin{aligned}
        \hat{X}_{t} & = x + \int_{0}^{t}b(\tau(s),\hat{X}_{\tau(s)},\hat{Y}_{\tau(s)})\dif s+ \int_{0}^{t} \langle \sigma(\tau(s),\hat{X}_{\tau(s)}), \dif W_{s}\rangle,          \\
        \hat{Y}_{t} & =\tilde{y}-\int_{0}^{t}f(\tau(s),\hat{X}_{\tau(s)},\hat{Y}_{\tau(s)},\tilde{Z}_{\tau(s)})\dif s+\int_{0}^{t} \langle \tilde{Z}_{\tau(s)},\dif W_{s}\rangle,
    \end{aligned}
    \right.
\end{equation}
with the aim of minimizing the objective functional
\[\hat{\J}(\tilde{y},\tilde{Z}_{.})=\left(\E|g(\hat{X}_{T})-\hat{Y}_{T}|^{2}\right)^{\frac{1}{2}}.\]
To characterize the error of time discretization, we define the modulus of continuity of the control \(\tilde{Z}_{.}\) as
\[\omega(\Delta t,\tilde{Z}_{.}):=\sup_{|t-s|\leq \Delta t,\,0\leq t,\,s\leq T}\left(\E|\tilde{Z}_{t}-\tilde{Z}_{s}|^{2}\right)^{\frac{1}{2}}.\]

 \begin{thm}
    \label{thm-dis-error}
    We assume Assumptions 1-2 hold.
    Let \((\tilde{X}_{t},\tilde{Y}_{t})\) and \((\hat{X}_{t},\hat{Y}_{t})\) be the state of system \eqref{eqn-fsde} and system \eqref{eqn-dis-fsde} respectively under the same control \((\tilde{y},\tilde{Z}_{.})\in \bR^{m}\times \cL^{2}_{\mathbf{F}}([0,T],\bR^{m\times d})\).
    We denote 
    \[|| (\tilde{y},\tilde{Z}_{.}) ||:=\left(\tilde{y}^{2}+ \int_{0}^{T}\E|\tilde{Z}_{t}|^{2}\dif t \right)^{\frac{1}{2}}.\]
    We have
    \[\sup_{0\leq t\leq T}\E|\tilde{X}_{t}-\hat{X}_{t}|+\sup_{0\leq t\leq T}\E| \tilde{Y}_{t}-\hat{Y}_{t}|\leq C\left[\omega(\Delta t,\tilde{Z}_{.})+\frac{1+\omega(\Delta t,\tilde{Z}_{.})+||(\tilde{y},\tilde{Z}_{.})||}{|\ln [\Delta t\vee \omega(\Delta t,\tilde{Z}_{.})]|}\right].\]
    Furthermore, we have
    \[\sup_{0\leq t\leq T}\E|\tilde{X}_{t}-\hat{X}_{t}|^{2}+\sup_{0\leq t\leq T}\E| \tilde{Y}_{t}-\hat{Y}_{t}|^{2}\leq C\left[\omega(\Delta t,\tilde{Z}_{.})+\omega(\Delta t,\tilde{Z}_{.})^{2}+\frac{1+\omega(\Delta t,\tilde{Z}_{.})+||(\tilde{y},\tilde{Z}_{.})||}{|\ln [\Delta t\vee \omega(\Delta t,\tilde{Z}_{.})]|}\right].\]
    Here the constant \(C\) does not depend on the choice of \((\tilde{y},\tilde{Z}_{.})\) and \(\Delta t\).
\end{thm}
It should be noted that the upper bound of the time discretization error is affected by the control used \((\tilde{y},\tilde{Z}_{.})\).
In practice, however we can limit the admissible control in a space of processes with sufficient time regularity.
To be more specific, for any constants \(\alpha,\,K>0\), we define the admissible control set as
\[\cA_{\alpha,K}:=\left\{ \xi_{t}\in \cL^{2}_{\mathbf{F}}([0,T],\bR^{m\times d})| \limsup_{\Delta t\to 0}\omega(\Delta t,\xi_{.})\Delta t^{-\alpha}<K \right\}.\]
By Theorems \ref{thm-control} and \ref{thm-dis-error}, we have the following discrete-time version posterior estimate.

\begin{thm}
    \label{thm-dis}
    We assume Assumptions 1-4 hold.
    Let \((X_{t},Y_{t})\) be the solution to FBSDE \eqref{eqn-fbsde} and \((\hat{X}_{t},\hat{Y}_{t})\) be the state of system \eqref{eqn-dis-fsde} under the control \((\tilde{y},\tilde{Z}_{.})\in\bR^{m}\times \cA_{\alpha,K}\).
    Then there exists constants \(C\), \(C_{\alpha,K}\), and \(N_{\alpha,K}\), such that for any \(\Delta t \leq T N_{\alpha,K}^{-1}\) we have
    \[\sup_{0\leq t\leq T}\E|X_{t}-\hat{X}_{t}|+\sup_{0\leq t\leq T}\E|Y_{t}-\hat{Y}_{t}|\leq C\hat{\J}(\tilde{y},\tilde{Z}_{.}) + C_{\alpha,K}\left[ |\ln \Delta t|^{-\frac{1}{2}}+\hat{\J}(\tilde{y},\tilde{Z}_{.}) |\ln \Delta t|^{-1}\right].\]
    Furthermore, we have
    \begin{align*}
        &\sup_{0\leq t\leq T}\E|X_{t}-\hat{X}_{t}|^{2}+\sup_{0\leq t\leq T}\E|Y_{t}-\hat{Y}_{t}|^{2}+\int_{0}^{T}\E|Z_{t}-\tilde{Z}_{t}|^{2}\dif t\nonumber\\
        \leq & C\left[\hat{\J}(\tilde{y},\tilde{Z}_{.})+\hat{\J}(\tilde{y},\tilde{Z}_{.})^{2}\right] +  C_{\alpha,K}\left[ |\ln \Delta t|^{-\frac{1}{2}}+\hat{\J}(\tilde{y},\tilde{Z}_{.}) |\ln \Delta t|^{-1}\right].
    \end{align*}
    Here the constants \(C\), \(C_{\alpha,K}\) and \(N_{\alpha,K}\) do not depend on the choice of \((\tilde{y},\tilde{Z}_{.})\) and \(C\) does not depend on \(\alpha\) and \(K\).
\end{thm} 

    \begin{rem}
        We point out that \(\cA_{\alpha,K}\) includes most cases of interest.
        Recall the solution to FBSDE \eqref{eqn-fbsde} satisfies \(Z_{t}=\nabla_{x}u(t,X_{t})\sigma(t,X_{t})\). 
        Then, we have \(Z_{.}\in\cA_{\alpha,K}\) for some \(K\) if \(\nabla_{x}u\sigma\) is \(\alpha\)-H\"older continuous.
        This condition does not exclude some important cases such as \(\sigma(x)=\sqrt{x}\).
        Moreover, given \(Z_{.}\in\cA_{\alpha,K}\), it is not hard to show
        \[\inf_{(\tilde{y},\tilde{Z}_{.})\in\bR \times \cA_{\alpha,K}}\hat{J}(\tilde{y},\tilde{Z}_{.})\leq C_{\alpha,K}|\ln\Delta t|^{-\frac{1}{2}}.\]
    \end{rem}

\subsection{ Deep BSDE Method}
Based on the Deep BSDE Method \cite{weinan2017deep}, we present our numerical algorithm.
For a given accuracy, Theorem \ref{thm-dis} enables us to approximate the solution to FBSDE \eqref{eqn-fbsde} by minimizing the objective functional \(\hat{\J}(\tilde{y},\tilde{Z}_{.})\).
The Euler-Maruyama scheme of the forward SDE \eqref{eqn-fsde} can be written as
\begin{equation}
    \left\{
    \begin{aligned}
        \hat{X}_{t_{k+1}}= & \hat{X}_{t_{k}}+b(t_{k},\hat{X}_{t_{k}},\hat{Y}_{t_{k}})\Delta t + \langle \sigma(t_{k},\hat{X}_{t_{k}}),\Delta W_{k}\rangle,       \\
        \hat{Y}_{t_{k+1}}= & \hat{Y}_{t_{k}}-f(t_{k},\hat{X}_{t_{k}},\hat{Y}_{t_{k}},\tilde{Z}_{t_{k}})\Delta t + \langle\tilde{Z}_{t_{k}}, \Delta W_{k}\rangle,
    \end{aligned}
    \right.
\end{equation}
where \(\Delta W_{k}\) denotes \(W_{t_{k+1}}-W_{t_{k}}\).

We represent the control \((\tilde{y},\tilde{Z}_{.})\) as \((\alpha,\Phi(\cdot,\hat{X}_{.};\beta))\) where \((\alpha,\beta)\) is the parameter to be optimized.
Here, the feedback control \(\Phi(t,x;\beta)\) is a forward neural network with uniformly bounded derivatives with respect to \(\beta\).
Therefore, the feedback control \(\tilde{Z}_{.}\) is contained in some \(\cA_{\alpha,K}\) throughout the training.
Thanks to the universal approximation capability of neural networks, we may approximate the optimal feedback control \(\nabla_{x}u(t,x)\sigma(t,x)\) by \(\Phi(t,x;\beta)\).
Our Deep BSDE network's forward propagation is written as
\begin{equation}
    \label{eqn-forward}
    \left\{
    \begin{aligned}
        \hat{X}_{t_{k+1}}= & \hat{X}_{t_{k}}+b(t_{k},\hat{X}_{t_{k}},\hat{Y}_{t_{k}})\Delta t + \langle \sigma(t_{k},\hat{X}_{t_{k}}), \Delta W_{k}\rangle, \\
        \hat{Y}_{t_{k+1}}= & \hat{Y}_{t_{k}}-f(t_{k},\hat{X}_{t_{k}},\hat{Y}_{t_{k}},\tilde{Z}_{t_{k}})\Delta t +\langle \tilde{Z}_{t_{k}}, \Delta W_{k}\rangle,       \\
        \tilde{Z}_{t_{k+1}}= & \Phi(t_{k+1},\hat{X}_{t_{k+1}};\beta),
    \end{aligned}
    \right.
\end{equation}
with initial value \(\hat{X}_{t_{0}}=x\), \(\hat{Y}_{t_{0}}=\alpha\), and \(\tilde{Z}_{t_{0}}=\Phi(t_{0},\hat{X}_{t_{0}};\beta)\).

We denote \(M\) as the batch size, i.e., we generate \(M\) paths of Brownian motion \(\left\{W_{t_{k}}^{i}\right\}_{1\leq i\leq M}\) in one batch.
For each path, we compute \((\hat{X}^{i}_{t_{k}},\hat{Y}^{i}_{t_{k}},\tilde{Z}^{i}_{t_{k}})\) iteratively using \eqref{eqn-forward}.
Then, we approximate the objective functional \(\E |g(\hat{X}_{T})-\hat{Y}_{T}|^{2}\) by Monte Carlo simulation.
The complete algorithm is shown below.

\begin{algorithm}[H]
    \DontPrintSemicolon
    \SetKwInput{KwInput}{Input}
    \SetKwInput{KwOutput}{Output}
    \KwInput{Initial parameters \((\alpha,\beta)\), learning rate \(\theta \);}
    \KwOutput{Optimal parameters \((\alpha^{\pi},\beta^{\pi})\);}
    \KwData{Brownian motion \(\left\{W_{t_{k}}^{i}\right\}_{1\leq i\leq M}\) samples.}
    \Repeat{the end condition.}{
    \(\hat{X}_{t_{0}}^{i}=x\), \(\hat{Y}_{t_{0}}^{i}=\alpha\), \(\tilde{Z}_{t_{0}}^{i}=\Phi(t_{0},\hat{X}_{t_{0}}^{i};\beta)\);\;
    \For{\(k=0\) \KwTo \(N-1\)} {
    \(\hat{X}_{t_{k+1}}^{i}=\hat{X}_{t_{k}}^{i}+b(t_{k},\hat{X}_{t_{k}}^{i},\hat{Y}_{t_{k}}^{i})\Delta t +\sigma(t_{k},\hat{X}_{t_{k}}^{i})\Delta W_{t_{k}}^{i}\);\;
    \(\hat{Y}_{t_{k+1}}^{i}=\hat{Y}_{t_{k}}^{i}-f(t_{k},\hat{X}_{t_{k}}^{i},\hat{Y}_{t_{k}}^{i},\tilde{Z}_{t_{k}}^{i})\Delta t+ \tilde{Z}_{t_{k}}^{i}\Delta W_{t_{k}}^{i}\);\;
    \(\tilde{Z}_{t_{k+1}}^{i}=\Phi(t_{k+1},\hat{X}_{t_{k+1}}^{i};\beta)\); \tcp*{forward propagation}
    }
    \((\alpha,\beta)=(\alpha,\beta)-\theta \nabla\frac{1}{M}\sum_{i=1}^{M}|g(\hat{X}_{T}^{i})-\hat{Y}_{T}^{i}|^{2}\);    \tcp*{parameters update}
    }
    \caption{The Deep BSDE Method}
\end{algorithm}

\section{Continuous-time Posterior Estimate}
This section contains the proof of Theorem \ref{thm-control}.
The proof was inspired by \cite{cvitanic2005steepest}.

From now on, we mainly consider a one-dimensional case that can be easily generalized to the multi-dimensional case.
For the simplicity of notations, we denote \(X_{t}-\tilde{X}_{t}\), \(Y_{t}-\tilde{Y}_{t}\), \(Z_{t}-\tilde{Z}_{t}\), \(b(t,X_{t},Y_{t})-b(t,\tilde{X}_{t},\tilde{Y}_{t})\), \(f(t,X_{t},Y_{t},Z_{t})-f(t,\tilde{X}_{t},\tilde{Y}_{t},\tilde{Z}_{t})\), \(\sigma(t,X_{t})-\sigma(t,\tilde{X}_{t})\) by \(\Delta X_{t}\), \(\Delta Y_{t}\), \(\Delta Z_{t}\), \(\Delta b_{t}\), \(\Delta f_{t}\), \(\Delta \sigma_{t}\)  respectively.
The constant \(C\) can be changed from one line to the next.

We present a series of Yamada-Watanabe type \cite{yamada1971uniqueness} functions \(\varphi_{m}\in C^{2}(\R)\) to approximate the absolute value function.
\(\varphi_{m}\) will be used repeatedly in our proof.
Let \(\varphi_{m}\) satisfy the following conditions:
\begin{itemize}
    \item \(\varphi_{m}(0)=\varphi_{m}'(0)=0,\)
    \item \(\varphi_{m}(x)=\varphi_{m}(-x),\)
    \item \(\text{supp}(\varphi_{m}'')\subseteq \left[ -\frac{2}{m},-\frac{1}{m^{2}} \right] \cup  \left[ \frac{1}{m^{2}},\frac{2}{m} \right],\)
    \item \(0\leq\varphi_{m}''(x)|x|\leq \frac{1}{\ln m},\)
    \item \(\int_{\frac{1}{m^{2}}}^{\frac{2}{m}}\varphi_{m}''(x)\dif x=1.\)
\end{itemize}
We can show that for any \(x\in\bR\) and \(m\in\bZ_{+}\)
\begin{itemize}
    \item \(|x|-\frac{2}{m}\leq \varphi_{m}(x) \leq |x|\),
    \item  \(-1\leq \varphi_{m}'(x) \leq 1\), 
    \item  \(0\leq \varphi_{m}''(x) \leq\frac{m^{2}}{\ln m}\).
\end{itemize}
Applying It\^{o}'s formula to \(|\tilde{Y}_{t}-u(t,\tilde{X}_{t})|^{2}\), we obtain
\begin{equation}
    \label{eqn-ito}
    \E |\tilde{Y}_{T}-u(T,\tilde{X}_{T})|^{2}=\E |\tilde{Y}_{t}-u(t,\tilde{X}_{t})|^{2} - 2\int_{t}^{T}\E \left[(\tilde{Y}_{s}-u(s,\tilde{X}_{s}))\alpha_{s} \right]\dif s + \int_{t}^{T}\E  \beta_{s}^{2}\dif s,
\end{equation}
where
\begin{equation*}
    \alpha_{s}=f(s,\tilde{X}_{s},\tilde{Y}_{s},\tilde{Z}_{s})+u_{t}(s,\tilde{X}_{s})+u_{x}(s,\tilde{X}_{s})b(s,\tilde{X}_{s},\tilde{Y}_{s})+\frac{1}{2}u_{xx}(s,\tilde{X}_{s})\sigma^{2}(s,\tilde{X}_{s})
\end{equation*}
and
\begin{equation*}
    \beta_{s}=\tilde{Z}_{s}-u_{x}(s,\tilde{X}_{s})\sigma(s,\tilde{X}_{s}).
\end{equation*}
Since \(u\) satisfies equation \eqref{eqn-decouple}, we have
\begin{align*}
    |\alpha_{s}| \leq & \left| f(s,\tilde{X}_{s},\tilde{Y}_{s},\tilde{Z}_{s})-f(s,\tilde{X}_{s},u(s,\tilde{X}_{s}),u_{x}(s,\tilde{X}_{s})\sigma(s,\tilde{X}_{s}))\right| \\
    &+                  \left| u_{x}(s,\tilde{X}_{s})b(s,\tilde{X}_{s},\tilde{Y}_{s}) -u_{x}(s,\tilde{X}_{s})b(s,\tilde{X}_{s},u(s,\tilde{X}_{s})) \right|               \\
    \leq              & \left( L_{f}+ L_{u}L_{b}\right) | \tilde{Y}_{s} - u(s,\tilde{X}_{s}) | +L_{f}\left|\beta_{s}\right|.
\end{align*}
Plugging the above estimate into \eqref{eqn-ito}, we obtain
\begin{align*}
    &\E|\tilde{Y}_{t}-u(t,\tilde{X}_{t})|^{2}\dif s + \int_{t}^{T}\E\beta^{2}_{s}\dif s \\
    =& \E |\tilde{Y}_{T}-u(T,\tilde{X}_{T})|^{2}+2\int_{t}^{T}\E \left[(\tilde{Y}_{s}-u(s,\tilde{X}_{s}))\alpha_{s} \right]\dif s\\
    \leq& \E|g(\tilde{X}_{T})-\tilde{Y}_{T}|^{2} + (1+L_{f})^{2}\int_{t}^{T}\E|\tilde{Y}_{s}-u(s,\tilde{X}_{s})|^{2}\dif s+ (1+L_{f})^{-2}\int_{t}^{T}\E\alpha_{s}^{2}\dif s\\
    \leq& \E|g(\tilde{X}_{T})-\tilde{Y}_{T}|^{2}+ \int_{t}^{T}\E\beta^{2}_{s}\dif s+C\int_{t}^{T}\E|\tilde{Y}_{s}-u(s,\tilde{X}_{s})|^{2}\dif s.
\end{align*}
Therefore, we have
\begin{align*}
    \E|\tilde{Y}_{t}-u(t,\tilde{X}_{t})|^{2}\dif s \leq \E|g(\tilde{X}_{T})-\tilde{Y}_{T}|^{2}+C\int_{t}^{T}\E|\tilde{Y}_{s}-u(s,\tilde{X}_{s})|^{2}\dif s.
\end{align*}
From Gronwall inequality, we show
\begin{equation}
    \label{ineq-Ty}
    \sup_{0\leq t \leq T} \E |\tilde{Y}_{t} - u(t,\tilde{X}_{t})|^{2} \leq C \E |g(\tilde{X}_{T})-\tilde{Y}_{T}|^{2}=C\J(\tilde{y},\tilde{Z}_{.})^{2}.
\end{equation}
Similarly, applying It\^{o}'s formula to \(\varphi_{m}(\Delta X_{t})\), we have
\[\E\varphi_{m}(\Delta X_{t})=\E\varphi_{m}(0)+\int_{0}^{t}\E\left[\varphi_{m}'(\Delta X_{s})\Delta b_{s} \right]\dif s + \frac{1}{2}\int_{0}^{t}\E\left[\varphi_{m}''(\Delta X_{s})\Delta \sigma^{2}_{s}\right]\dif s.\]
Notice that
\begin{align}
    \label{est-b}
    |\Delta b_{s}|= & |b(s,X_{s},Y_{s})-b(s,\tilde{X}_{s},\tilde{Y}_{s})|                                     \nonumber               \\
    \leq            & L_{b}\left( |X_{s}-\tilde{X}_{s}|+|Y_{s}-u(s,\tilde{X}_{s})|+|\tilde{Y}_{s}-u(s,\tilde{X}_{s})|\right) \nonumber \\
    \leq            & (L_{b}+L_{u}L_{b})|X_{s}-\tilde{X}_{s}| +L_{b}|\tilde{Y}_{s}-u(s,\tilde{X}_{s})|
\end{align}
and 
\begin{align}
    \label{est-sig}
    |\Delta\sigma_{s}|^{2}\leq L_{\sigma}|X_{s}-\tilde{X}_{s}|.
\end{align}
Then, combining the properties of \(\varphi_{m}\) we get
\begin{align*}
    \E |\Delta X_{t}|\leq & \frac{2}{m}+\E \varphi_{m}(\Delta X_{t})                                                                                          \\
    \leq                  & \frac{2}{m}+ \frac{C}{\ln m} + C\sup_{0\leq t\leq T}\E|\tilde{Y}_{t}-u(t,\tilde{X}_{t})| + C\int_{0}^{t}\E|\Delta X_{s}| \dif s.
\end{align*}
Here the constant \(C\) does not depend on \(m\).
Let \(m\to \infty\),
it follows from Gronwall inequality that
\begin{equation}
    \label{ineq-Tx}
    \sup_{0\leq t\leq T}\E|\Delta X_{t}|\leq C \sup_{0\leq t\leq T}\E|\tilde{Y}_{t}-u(t,\tilde{X}_{t})|.
\end{equation}
In view of estimate \eqref{ineq-Ty}, we obtain
\begin{equation}
    \label{ineq-x1}
    \sup_{0\leq t\leq T}\E|\Delta X_{t}|\leq C \J(\tilde{y},\tilde{Z}_{t}).
\end{equation}
Again, by Assumptions 3, 4, and estimate \eqref{ineq-Ty}, we have
\begin{align*}
    \sup_{0\leq t\leq T}\E|\Delta Y_{t}|\leq & \sup_{0\leq t\leq T} \E|u(t,X_{t})-u(t,\tilde{X}_{t})|+\sup_{0\leq t\leq T}\E|\tilde{Y}_{t}-u(t,\tilde{X}_{t})|        \\
    \leq                                     & L_{u}\sup_{0\leq t\leq T} \E |\Delta X_{t}| + \sqrt{\sup_{0\leq t\leq T}\E|\tilde{Y}_{t}-u(t,\tilde{X}_{t})|^{2}} \\
    \leq                                     & C\J(\tilde{y},\tilde{Z}_{t}).
\end{align*}

Now, we begin to estimate the second moment of errors.
Applying It\^{o}'s formula to \(|\Delta X_{t}|^{2}\) we obtain
\begin{align*}
    \E |\Delta X_{t}|^{2}= &2\int_{0}^{t}\E\left[ \Delta X_{s}\Delta b_{s} \right]\dif s +\int_{0}^{t}\E |\Delta \sigma_{s}|^{2}\dif s\\
    \leq & \int_{0}^{t}\E|\Delta X_{s}|^{2}\dif s + \int_{0}^{t}\E|\Delta b_{s}|^{2}\dif s+\int_{0}^{t}\E |\Delta \sigma_{s}|^{2}\dif s.
\end{align*}
Plugging estimates \eqref{ineq-Ty}, \eqref{est-b}, \eqref{est-sig}, and \eqref{ineq-x1}, we have
\begin{equation*}
    \E |\Delta X_{t}|^{2}\leq C\int_{0}^{t}\E|\Delta X_{s}|^{2}\dif s + C\left[\J(\tilde{y},\tilde{Z}_{t})+\J(\tilde{y},\tilde{Z}_{t})^{2}\right].
\end{equation*}
By Gronwall inequality, we show
\begin{equation}
    \label{ineq-x2}
    \sup_{0\leq t\leq T}\E |\Delta X_{t}|^{2}\leq C\left[\J(\tilde{y},\tilde{Z}_{t})+\J(\tilde{y},\tilde{Z}_{t})^{2}\right].
\end{equation}
Moreover, by Assumption 4 and estimate \eqref{ineq-Ty}, we have 
\begin{align}
    \label{ineq-y2}
    \sup_{0\leq t\leq T}\E|\Delta Y_{t}|^{2}\leq & 2\sup_{0\leq t\leq T}\E|Y_{t}-u(t,\tilde{X}_{t})|^{2}+2\sup_{0\leq t\leq T}\E |\tilde{Y}_{t}-u(t,\tilde{X}_{t})|^{2}\nonumber\\
    \leq &2L_{u}^{2}\sup_{0\leq t\leq T}\E|\Delta X_{t}|^{2}+C \E|g(\tilde{X}_{T})-\tilde{Y}_{T}|^{2} \nonumber\\
    \leq &C\left[\J(\tilde{y},\tilde{Z}_{t})+\J(\tilde{y},\tilde{Z}_{t})^{2}\right].
\end{align}
Applying It\^{o}'s formula to \(|\Delta Y_{t}|^{2}\) we get
\begin{equation*}
    \E |\Delta Y_{T}|^{2} = \E |\Delta Y_{0}|^{2} -2\int_{0}^{T}\E\left[ \Delta f_{s}\Delta Y_{s} \right]\dif s + \int_{0}^{T} \E|\Delta Z_{s}|^{2} \dif s.
\end{equation*}
Combining with the fact that 
\begin{equation*}
    |\Delta f_{s}|\leq L_{f}\left( |\Delta X_{s}|+|\Delta Y_{s}|+ |\Delta Z_{s}| \right),
\end{equation*}
we show 
\begin{align*}
    \int_{0}^{T} \E|\Delta Z_{s}|^{2} \dif s \leq& 2\int_{0}^{T}\E\left[ \Delta f_{s}\Delta Y_{s} \right]\dif s + \sup_{0\leq t\leq T}\E|\Delta Y_{t}|^{2}\\
    \leq& (2L_{f}+1)^{2}\int_{0}^{T}\E|\Delta Y_{s}|^{2}\dif s + (2L_{f}+1)^{-2}\int_{0}^{T}\E|\Delta f_{s}|^{2}\dif s+ \sup_{0\leq t\leq T}\E|\Delta Y_{t}|^{2}\\
    \leq& C\sup_{0\leq t\leq T}\E|\Delta X_{t}|^{2} +C\sup_{0\leq t\leq T}\E|\Delta Y_{t}|^{2}+\frac{1}{2}\int_{0}^{T} \E|\Delta Z_{s}|^{2} \dif s.
\end{align*}
Together with estimates \eqref{ineq-x2} and \eqref{ineq-y2}, we have
\begin{equation}
    \int_{0 }^{T}\E |\Delta Z_{t}|^{2}\dif s\leq C\left[\J(\tilde{y},\tilde{Z}_{t})+\J(\tilde{y},\tilde{Z}_{t})^{2}\right].
\end{equation}

\section{Time Discretization Error}
This section proves Theorem \ref{thm-dis-error}.

Recall that \((\hat{X}_{t},\hat{Y}_{t})\) satisfies the following forward SDE
\begin{equation}
    \label{eqn-sys}
    \left\{
    \begin{aligned}
        \hat{X}_{t} & = x + \int_{0}^{t}b(\tau(s),\hat{X}_{\tau(s)},\hat{Y}_{\tau(s)})\dif s+ \int_{0}^{t}\sigma(\tau(s),\hat{X}_{\tau(s)})\dif W_{s},           \\
        \hat{Y}_{t} & =\tilde{y}-\int_{0}^{t}f(\tau(s),\hat{X}_{\tau(s)},\hat{Y}_{\tau(s)},\tilde{Z}_{\tau(s)})\dif s+\int_{0}^{t}\tilde{Z}_{\tau(s)}\dif W_{s},
    \end{aligned}
    \right.
\end{equation}
and
\[\Delta t=\frac{T}{N},\, t_{k}=\frac{kT}{N},\,\Delta W_{k}=W_{t_{k+1}}-W_{t_{k}}.\]
We have
\begin{equation}
    \label{eqn-step}
    \left\{
    \begin{aligned}
        \hat{X}_{t_{k+1}}= & \hat{X}_{t_{k}}+b(t_{k},\hat{X}_{t_{k}},\hat{Y}_{t_{k}})\Delta t + \sigma(t_{k},\hat{X}_{t_{k}})\Delta W_{k},      \\
        \hat{Y}_{t_{k+1}}= & \hat{Y}_{t_{k}}-f(t_{k},\hat{X}_{t_{k}},\hat{Y}_{t_{k}},\tilde{Z}_{t_{k}})\Delta t +\tilde{Z}_{t_{k}}\Delta W_{k}.
    \end{aligned}
    \right.
\end{equation}
By \(b_{0}\), \(\sigma_{0}\), and \(f_{0}\) we denote \(b(0,0,0)\), \(\sigma(0,0)\), and \(f(0,0,0,0)\), respectively.

\begin{lem}
    \label{lem-bound}
    We assume Assumption 2 holds.
    Let \((\hat{X}_{t},\hat{Y}_{t})\) be the state of system \eqref{eqn-sys} under the control \((\tilde{y},\tilde{Z}_{.})\in \bR^{m}\times \cL^{2}_{\mathbf{F}}([0,T],\bR^{m\times d})\).
    We have for some constant \(C\)
    \[\sup_{0\leq k\leq N}\E|\hat{X}_{t_{k}}|^{2}+\sup_{0\leq k\leq N}\E|\hat{Y}_{t_{k}}|^{2}\leq C\left[1+\omega(\Delta t,\tilde{Z}_{.})^{2}+||(\tilde{y},\tilde{Z}_{.})||^{2}\right].\]
    Here the constant \(C\) does not depend on the choice of \((\tilde{y},\tilde{Z}_{.})\) and \(N\).
\end{lem}

\begin{proof}
    We square both sides of the equation \eqref{eqn-step} and take expectations.
    For \(\hat{X}_{t_{k}}\), we have
        \begin{align}
            \label{ineq-step-x}
            \E|\hat{X}_{t_{k+1}}|^{2}\leq & (1+\Delta t)\E \hat{X}_{t_{k}}^{2} + (\Delta t +\Delta t^{2})\E \left[b^{2}(t_{k},\hat{X}_{t_{k}},\hat{Y}_{t_{k}})\right] + \E\left[\sigma^{2}(t_{k},\hat{X}_{t_{k}})\Delta W_{k}^{2}\right] \nonumber\\
            \leq                        & (1+\Delta t)\E |\hat{X}_{t_{k}}|^{2} + C\Delta t \E \left[b^{2}_{0}+L_{b}^{2}T^{2}+L_{b}^{2}|\hat{X}_{t_{k}}|^{2}+L_{b}^{2}|\hat{Y}_{t_{k}}|^{2}\right]               \nonumber\\
            &+                            C\Delta t\E\left[\sigma^{2}_{0}+L_{\sigma}T+L_{\sigma}|\hat{X}_{t_{k}}|\right]                                                                                 \nonumber \\
            \leq                        & (1+C\Delta t)\E|\hat{X}_{t_{k}}|^{2} + C\Delta t \E|\hat{Y}_{t_{k}}|^{2} + C\Delta t +C\Delta t \left( \E|\hat{X}_{t_{k}}|^{2} \right)^{\frac{1}{2}}\nonumber\\
            \leq &(1+C\Delta t)\E|\hat{X}_{t_{k}}|^{2} + C\Delta t \E|\hat{Y}_{t_{k}}|^{2} + C\Delta t.
        \end{align}
    For \(\hat{Y}_{t_{k}}\), we have
        \begin{align}
            \label{ineq-step-y}
            &\E |\hat{Y}_{t_{k+1}}|^{2}\nonumber\\
            \leq & (1+\Delta t)\E|\hat{Y}_{t_{k}}|^{2}+  (\Delta t +\Delta t^{2})\E \left[f^{2}(t_{k},\hat{X}_{t_{k}},\hat{Y}_{t_{k}},\tilde{Z}_{t_{k}})\right] +\E \left[|\tilde{Z}_{t_{k}}|^{2}\Delta W_{k}^{2} \right]                             \nonumber               \\
            \leq                       & (1+\Delta t)\E |\hat{Y}_{t_{k}}|^{2} + C\Delta t \E \left[f^{2}_{0}+L_{f}^{2}T^{2}+L_{f}^{2}|\hat{X}_{t_{k}}|^{2}+L_{f}^{2}|\hat{Y}_{t_{k}}|^{2}+L_{f}^{2}|\tilde{Z}_{t_{k}}|^{2}\right]+\Delta t\E |\tilde{Z}_{t_{k}}|^{2}  \nonumber\\
            \leq                       & (1+C\Delta t)\E|\hat{Y}_{t_{k}}|^{2} + C\Delta t \E|\hat{X}_{t_{k}}|^{2} + C\Delta t+C\Delta t\E|\tilde{Z}_{t_{k}}|^{2}\nonumber\\
            \leq & (1+C\Delta t)\E|\hat{Y}_{t_{k}}|^{2} + C\Delta t \E|\hat{X}_{t_{k}}|^{2} + C\Delta t+ C\int_{t_{k}}^{t_{k+1}}\E|\tilde{Z}_{t}|^{2}\dif t +C\Delta t \omega(\Delta t,\tilde{Z}_{.})^{2}.
        \end{align}
    Combined with estimate \eqref{ineq-step-x}, we have
    \begin{equation*}
        \E \left[|\hat{X}_{t_{k+1}}|^{2}+|\hat{Y}_{t_{k+1}}|^{2}\right]\leq (1+C\Delta t)\E \left[|\hat{X}_{t_{k}}|^{2}+|\hat{Y}_{t_{k}}|^{2}\right] +C\Delta t+ C\int_{t_{k}}^{t_{k+1}}\E|\tilde{Z}_{t}|^{2}\dif t +C\Delta t \omega(\Delta t,\tilde{Z}_{.})^{2}.
    \end{equation*}
    Here, we define an increasing sequence of numbers \(\left\{ a_{k}\right\}\) satisfying
    \begin{equation*}
        \left\{
        \begin{aligned}
            a_{k+1}= & (1+C\Delta t)a_{k}+C\left[\Delta t+\int_{t_{k}}^{t_{k+1}}\E|\tilde{Z}_{t}|^{2}\dif t +\Delta t \omega(\Delta t,\tilde{Z}_{.})^{2}\right],      \\
            a_{0}= & x^{2}+\tilde{y}^{2}.
        \end{aligned}
        \right.
    \end{equation*}
    Thus, we have
    \begin{align*}
        a_{N}\leq & (1+C\Delta t)^{N}a_{0}+ C\sum_{k=0}^{N-1}(1+C\Delta t )^{k}\left[\Delta t+\int_{t_{k}}^{t_{k+1}}\E|\tilde{Z}_{t}|^{2}\dif t +\Delta t \omega(\Delta t,\tilde{Z}_{.})^{2}\right]\\
             \leq & C\left[1+\omega(\Delta t,\tilde{Z}_{.})^{2}+\tilde{y}^{2}+\int_{0}^{T}\E |\tilde{Z}_{t}|^{2}\dif t\right]\\
             =&C\left[1+\omega(\Delta t,\tilde{Z}_{.})^{2}+||(\tilde{y},\tilde{Z}_{.})||^{2}\right].
    \end{align*}
    We conclude by \[\E\left[ |\hat{X}_{t_{k}}|^{2}+|\hat{Y}_{t_{k}}|^{2} \right]\leq a_{k}\leq a_{N}.\]
\end{proof}

\begin{lem}
    \label{lem-continuous}
    We assume Assumption 2 holds.
    Let \((\hat{X}_{t},\hat{Y}_{t})\) be the state of system \eqref{eqn-sys} under the control \((\tilde{y},\tilde{Z}_{.})\in \bR^{m}\times \cL^{2}_{\mathbf{F}}([0,T],\bR^{m\times d})\).
    We have for some constant \(C\)
    \[\sup_{0\leq t\leq T}\E|\hat{X}_{t}-\hat{X}_{\tau(t)}|^{2}+\int_{0}^{T}\E|\hat{Y}_{t}-\hat{Y}_{\tau(t)}|^{2}\dif t\leq C \Delta t\left[1+\omega(\Delta t,\tilde{Z}_{.})^{2}+||(\tilde{y},\tilde{Z}_{.})||^{2}\right].\]
    Here the constant \(C\) does not depend on the choice of \((\tilde{y},\tilde{Z}_{.})\) and \(N\).
\end{lem}

\begin{proof}
    Without loss of generality, we assume \(t_{k}\leq t\leq t_{k+1}\).
    It follows from Lemma \ref{lem-bound} and Assumption 2 that
    \begin{align*}
         \E|\hat{X}_{t}-\hat{X}_{t_{k}}|^{2}\leq & C\Delta t \int_{t_{k}}^{t}\E\left[b^{2}(\tau(s),\hat{X}_{\tau(s)},\hat{Y}_{\tau(s)})\right]\dif s +C\E \left|\int_{t_{k}}^{t}\sigma(\tau(s),\hat{X}_{\tau(s)}) \dif W_{s}\right|^{2} \\
        \leq                                                   & C\Delta t \int_{t_{k}}^{t}\E\left[b^{2}(\tau(s),\hat{X}_{\tau(s)},\hat{Y}_{\tau(s)})\right]\dif s+C\int_{t_{k}}^{t} \E\left[\sigma^{2}(\tau(s),\hat{X}_{\tau(s)})\right]\dif s                              \\
        \leq                                                   & C\Delta t \int_{t_{k}}^{t}\E \left[|b_{0}|+L_{b}T+L_{b}|\hat{X}_{\tau(s)}|+L_{b}|\hat{Y}_{\tau(s)}|\right]^{2}\dif s                                                                                 \\
        &+                                                       C\int_{t_{k}}^{t}\E\left[\sigma^{2}_{0}+L_{\sigma}T+L_{\sigma}|\hat{X}_{\tau(s)}|\right]\dif s                                                                                                                 \\
        \leq                                                   & C\Delta t\left[1+\omega(\Delta t,\tilde{Z}_{.})^{2}+||(\tilde{y},\tilde{Z}_{.})||^{2}\right].
    \end{align*}
    Similarly,
    \begin{align*}
         \E|\hat{Y}_{t}-\hat{Y}_{t_{k}}|^{2}\leq & C\Delta t\int_{t_{k}}^{t}\E\left[f^{2}(\tau(s),\hat{X}_{\tau(s)},\hat{Y}_{\tau(s)},\tilde{Z}_{\tau(s)})\right]\dif s +C \E \left|\int_{t_{k}}^{t}\tilde{Z}_{\tau(s)} \dif W_{s}\right|^{2}\\
        \leq                                                   & C\Delta t\int_{t_{k}}^{t}\E\left[f^{2}(\tau(s),\hat{X}_{\tau(s)},\hat{Y}_{\tau(s)},\tilde{Z}_{\tau(s)})\right]\dif s+C\int_{t_{k}}^{t} \E|\tilde{Z}_{\tau(s)}|^{2}\dif s                              \\
        \leq                                                   & C \Delta t\int_{t_{k}}^{t}\E \left[|f_{0}|+L_{f}T+L_{f}|\hat{X}_{\tau(s)}|+L_{y}|\hat{Y}_{\tau(s)}|+L_{f}|\tilde{Z}_{\tau(s)}|\right]^{2}\dif s +C\Delta t \E|\tilde{Z}_{t_{k}}|^{2}                                            \\
        \leq                                                   & C\Delta t\left[1+\omega(\Delta t,\tilde{Z}_{.})^{2}+||(\tilde{y},\tilde{Z}_{.})||^{2}\right]+C\Delta t \E|\tilde{Z}_{t_{k}}|^{2}.
    \end{align*}
    Therefore, we obtain 
    \begin{align*}
        \int_{0}^{T}\E |\hat{Y}_{t}-\hat{Y}_{\tau(t)}|^{2}\dif t\leq &C\Delta t\left[1+\omega(\Delta t,\tilde{Z}_{.})^{2}+||(\tilde{y},\tilde{Z}_{.})||^{2}\right] +C\Delta t \int_{0}^{T}\E|\tilde{Z}_{\tau(t)}|^{2}\dif t\\
        \leq &C\Delta t\left[1+\omega(\Delta t,\tilde{Z}_{.})^{2}+||(\tilde{y},\tilde{Z}_{.})||^{2}\right].
    \end{align*}
\end{proof}

We denote \(\tilde{X}_{t}-\hat{X}_{t}\), \(\tilde{Y}_{t}-\hat{Y}_{t}\), \(b(t,\tilde{X}_{t},\tilde{Y}_{t})-b(\tau(t),\hat{X}_{\tau(t)},\hat{Y}_{\tau(t)})\),  \(\sigma(t,\tilde{X}_{t})-\sigma(\tau(t),\hat{X}_{\tau(t)})\), \(f(t,\tilde{X}_{t},\tilde{Y}_{t},\tilde{Z}_{t})-f(\tau(t),\hat{X}_{\tau(t)},\hat{Y}_{\tau(t)},\tilde{Z}_{\tau(t)})\) by \(\hat{\Delta} X_{t}\), \(\hat{\Delta} Y_{t}\), \(\hat{\Delta} b_{t}\), \(\hat{\Delta} \sigma_{t}\), \(\hat{\Delta} f_{t}\), respectively.
From It\^{o}'s formula, we have
\[\E\varphi_{m}(\hat{\Delta} X_{t})=\E\varphi_{m}(0)-\int_{0}^{t}\E\left[\varphi_{m}'(\hat{\Delta} X_{s})\hat{\Delta} b_{s} \right]\dif s + \frac{1}{2}\int_{0}^{t}\E\left[\varphi_{m}''(\hat{\Delta} X_{s})\hat{\Delta} \sigma^{2}_{s}\right]\dif s.\]
For the second term, by Lemma \ref{lem-continuous}, we have
\begin{align}
    \label{est-x1}
         & \left| \int_{0}^{t}\E\left[\varphi_{m}'(\hat{\Delta} X_{s})\hat{\Delta} b_{s} \right]\dif s  \right|       \nonumber                                                                                         \\
    \leq & \int_{0}^{t}\left\{\E| b(s,\tilde{X}_{s},\tilde{Y}_{s})-b(s,\hat{X}_{s},\hat{Y}_{s})|+\E|b(s,\hat{X}_{s},\hat{Y}_{s})-b(\tau(s),\hat{X}_{\tau(s)},\hat{Y}_{\tau(s)})|\right\}\dif s \nonumber\\
    \leq & L_{b}\int_{0}^{t}\E \left[|\hat{\Delta} X_{s}|+|\hat{\Delta} Y_{s}|+|s-\tau(s)|+|\hat{X}_{s}-\hat{X}_{\tau(s)}|+|\hat{Y}_{s}-\hat{Y}_{\tau(s)}|\right]\dif s                  \nonumber   \\
    \leq & C\int_{0}^{t}\E|\hat{\Delta} X_{s}|\dif s + C\int_{0}^{t}\E|\hat{\Delta} Y_{s}|\dif s + C\Delta t +C \sup_{0\leq t\leq T}\E|\hat{X}_{t}-\hat{X}_{\tau(t)}| +C \int_{0}^{T}\E |\hat{Y}_{t}-\hat{Y}_{\tau(t)}|\dif t \nonumber\\
    \leq& C\int_{0}^{t}\E|\hat{\Delta} X_{s}|\dif s + C\int_{0}^{t}\E|\hat{\Delta} Y_{s}|\dif s+C\Delta t^{\frac{1}{2}}\left[1+\omega(\Delta t,\tilde{Z}_{.})+||(\tilde{y},\tilde{Z}_{.})||\right].
\end{align}
For the third term, it also follows from Lemma \ref{lem-continuous} that
\begin{align}
    \label{est-x2}
         & \left|\frac{1}{2}\int_{0}^{t}\E\left[\varphi_{m}''(\hat{\Delta} X_{s})\hat{\Delta} \sigma^{2}_{s}\right]\dif s \right|                   \nonumber                                                                                                             \\
    \leq & \int_{0}^{t}\left\{\E\left[\varphi_{m}''(\hat{\Delta} X_{s})|\sigma(s,\tilde{X}_{s})-\sigma(s,\hat{X}_{s}) |^{2}\right]+\E\left[\varphi_{m}''(\hat{\Delta} X_{s})|\sigma(s,\hat{X}_{s})-\sigma(\tau(s),\hat{X}_{\tau(s)})|^{2}\right] \right\}\dif s\nonumber \\
    \leq & L_{\sigma}\int_{0}^{t}\E \left[\varphi_{m}''(\hat{\Delta} X_{s})|\hat{\Delta} X_{s}|\right]\dif s + L_{\sigma}\left\|\varphi_{m}''\right\|_{\infty}\int_{0}^{t}\left[|s-\tau(s)|+\E|\hat{X}_{s}-\hat{X}_{\tau(s)}|\right]\dif s                     \nonumber           \\
    \leq & \frac{C}{\ln m} +\frac{Cm^{2}\Delta t^{\frac{1}{2}}}{\ln m}\left[1+\omega(\Delta t,\tilde{Z}_{.})+||(\tilde{y},\tilde{Z}_{.})||\right].
\end{align}
Applying It\^{o}'s formula to \(\varphi_{m}(\hat{\Delta} Y_{t})\), we get
\[\E\varphi_{m}(\hat{\Delta} Y_{t})=\E\varphi_{m}(0)+\int_{0}^{t}\E\left[\varphi_{m}'(\hat{\Delta} Y_{s})\hat{\Delta} f_{s}\right] \dif s + \frac{1}{2}\int_{0}^{t}\E\left[\varphi_{m}''(\hat{\Delta} Y_{s})| \tilde{Z}_{s}-\tilde{Z}_{\tau(s)} |^{2}\right]\dif s.\]
For the second term, we have
\begin{align} 
    \label{est-y1}
    &\left| \int_{0}^{t}\E\left[\varphi_{m}'(\hat{\Delta} Y_{s})\hat{\Delta} f_{s} \right]\dif s  \right|\nonumber\\
    \leq & \int_{0}^{t}\left\{ \E|f(s,\tilde{X}_{s},\tilde{Y}_{s},\tilde{Z}_{s})-f(s,\hat{X}_{s},\hat{Y}_{s},\tilde{Z}_{s})|+\E|f(s,\hat{X}_{s},\hat{Y}_{s},\tilde{Z}_{s})-f(\tau(s),\hat{X}_{\tau(s)},\hat{Y}_{\tau(s)},\tilde{Z}_{\tau(s)})| \right\}\dif s\nonumber\\
    \leq & L_{f}\int_{0}^{t}\E\left[ |\hat{\Delta }X_{s}|+|\hat{\Delta}Y_{s}|+|s-\tau(s)|+|\hat{X}_{s}-\hat{X}_{\tau_{s}}|+|\hat{Y}_{s}-\hat{Y}_{\tau(s)}|+|\tilde{Z}_{s}-\tilde{Z}_{\tau(s)}| \right]\dif s\nonumber\\
    \leq & C\int_{0}^{t}\E|\hat{\Delta} X_{s}|\dif s + C\int_{0}^{t}\E|\hat{\Delta} Y_{s}|\dif s+C\Delta t^{\frac{1}{2}}\left[1+\omega(\Delta t,\tilde{Z}_{.})+||(\tilde{y},\tilde{Z}_{.})||\right] +C\omega(\Delta t,\tilde{Z}_{.}).
\end{align}
For the third term,
\begin{align}
    \label{est-y2}
    \left|\frac{1}{2}\int_{0}^{t}\E\left[\varphi_{m}''(\hat{\Delta} Y_{s})| \tilde{Z}_{s}-\tilde{Z}_{\tau(s)}|^{2}\right]\dif s \right|\leq C\left\|\varphi_{m}''\right\|_{\infty}\omega(\Delta t,\tilde{Z}_{.})^{2}\leq C\frac{m^{2}}{\ln m}\omega(\Delta t,\tilde{Z}_{.})^{2}.
\end{align}
Notice the fact that \[|x|\leq\varphi_{m}(x)+\frac{2}{m}.\]
Hence, by estimates \eqref{est-x1}, \eqref{est-x2}, \eqref{est-y1}, and \eqref{est-y2}, we have
\begin{align*}
    &\E|\hat{\Delta} X_{t}|+\E|\hat{\Delta} Y_{t}|\\
    \leq & \frac{4}{m}+\E\varphi_{m}(\hat{\Delta} X_{t})+\E\varphi_{m}(\hat{\Delta} Y_{t})                                                                         \\
    \leq                                  & \frac{C}{\ln m}+ \frac{Cm^{2}\Delta t^{\frac{1}{2}}}{\ln m}\left[1+\omega(\Delta t,\tilde{Z}_{.})+||(\tilde{y},\tilde{Z}_{.})||\right]+C\omega(\Delta t,\tilde{Z}_{.})+C\frac{m^{2}}{\ln m}\omega(\Delta t,\tilde{Z}_{.})^{2}\\
    &+C\int_{0}^{t}\E|\hat{\Delta} X_{s}|+\E|\hat{\Delta} Y_{s}|\dif s.
\end{align*}
We set \(m=\lfloor \Delta t^{-\frac{1}{4}} \wedge \omega(\Delta t,\tilde{Z}_{.})^{-1}\rfloor\).
From Gronwall inequality, we have
\[\sup_{0\leq t\leq T}\E|\hat{\Delta} X_{t}|+\sup_{0\leq t\leq T}\E|\hat{\Delta} Y_{t}|\leq C\left[\omega(\Delta t,\tilde{Z}_{.})+\frac{1+\omega(\Delta t,\tilde{Z}_{.})+||(\tilde{y},\tilde{Z}_{.})||}{|\ln [\Delta t\vee \omega(\Delta t,\tilde{Z}_{.})]|}\right].\]
Similarly, we apply It\^{o}'s formula to \(|\hat{\Delta} X_{t}|^{2}\) and \(|\hat{\Delta} Y_{t}|^{2}\).
Following standard arguments, we have
\[\sup_{0\leq t\leq T}\E|\hat{\Delta} X_{t}|^{2}+\sup_{0\leq t\leq T}\E|\hat{\Delta} Y_{t}|^{2}\leq C\left[\omega(\Delta t,\tilde{Z}_{.})+\omega(\Delta t,\tilde{Z}_{.})^{2}+\frac{1+\omega(\Delta t,\tilde{Z}_{.})+||(\tilde{y},\tilde{Z}_{.})||}{|\ln [\Delta t\vee \omega(\Delta t,\tilde{Z}_{.})]|}\right].\]
Hence, we prove Theorem \ref{thm-dis-error}.

\section{Discrete-time Posterior Estimate}
    This section proves Theorem \ref{thm-dis}.
    Recall that we assume
    \[\tilde{Z}_{.}\in\cA_{\alpha,K}=\left\{ \xi_{t}\in \cL^{2}_{\mathbf{F}}([0,T],\bR^{m\times d})| \limsup_{\Delta t\to 0}\omega(\Delta t,\xi_{.})\Delta t^{-\alpha}<K \right\}.\]
    Without loss of generality, we may assume \(N_{\alpha,K}\) is sufficiently large such that
    \[\omega(\Delta t,\tilde{Z}_{.})\leq K \Delta t^{\alpha}\leq |\ln\Delta t|^{-1}<1,\]
    for some \(K>0\) and any \(\Delta t \leq TN_{\alpha,K}^{-1}\).
    By Assumption 2, we have 
    \begin{align*}
        \J(\tilde{y},\tilde{Z}_{.})^{2}\leq C\hat{\J}(\tilde{y},\tilde{Z}_{.})^{2}+C\E|X_{T}-\hat{X}_{T}|^{2}+C\E|Y_{T}-\hat{Y}_{T}|^{2}.
    \end{align*}
    By Theorem \ref{thm-dis-error}, we have
    \begin{align}
        \label{est-cor1}
        &\sup_{0\leq t\leq T}\E| X_{t}-\hat{X}_{t}|^{2}+\sup_{0\leq t\leq T}\E|Y_{t}-\hat{Y}_{t}|^{2}\nonumber\\
        \leq& C\left[\omega(\Delta t,\tilde{Z}_{.})+\omega(\Delta t,\tilde{Z}_{.})^{2}+\frac{1+\omega(\Delta t,\tilde{Z}_{.})+||(\tilde{y},\tilde{Z}_{.})||}{|\ln [\Delta t\vee \omega(\Delta t,\tilde{Z}_{.})]|}\right]\nonumber\\
        \leq& C_{\alpha,K}\left[ 1+||(\tilde{y},\tilde{Z}_{.})|| \right]|\ln \Delta t|^{-1}.
    \end{align}
    Similarly,
    \begin{align}
        \sup_{0\leq t\leq T}\E| X_{t}-\hat{X}_{t}|+\sup_{0\leq t\leq T}\E|Y_{t}-\hat{Y}_{t}|\leq C_{\alpha,K}\left[ 1+||(\tilde{y},\tilde{Z}_{.})|| \right]|\ln \Delta t|^{-1}.
    \end{align}
    By Theorem \ref{thm-control}, we obtain
    \begin{align}
        \label{est-cor2}
        ||(\tilde{y},\tilde{Z}_{.})||^{2}\leq &C \left[|Y_{0}-\tilde{Y}_{0}|^{2}+\int_{0}^{T}\E|Z_{t}-\tilde{Z}_{t}|^{2}\dif t\right]+C \left[|Y_{0}|^{2}+\int_{0}^{T}\E |Z_{t}|^{2}\dif t\right]\nonumber\\
        \leq& C\left[ 1+\J(\tilde{y},\tilde{Z}_{.})^{2} \right].
    \end{align}
    Combining above estimates and Assumption 2, we have 
    \begin{align*}
        \J(\tilde{y},\tilde{Z}_{.})^{2}\leq& C\hat{\J}(\tilde{y},\tilde{Z}_{.})^{2}+C\E|X_{T}-\hat{X}_{T}|^{2}+C\E|Y_{T}-\hat{Y}_{T}|^{2}\\
        \leq &C\hat{\J}(\tilde{y},\tilde{Z}_{.})^{2} +C_{\alpha,K}\left[ 1+||(\tilde{y},\tilde{Z}_{.})|| \right]|\ln \Delta t|^{-1}\\
        \leq &C\hat{\J}(\tilde{y},\tilde{Z}_{.})^{2} +C_{\alpha,K}\left[ 1 + \J(\tilde{y},\tilde{Z}_{.})\right]|\ln \Delta t|^{-1}.
    \end{align*}
    This implies
    \begin{equation} \J(\tilde{y},\tilde{Z}_{.})\leq C\hat{\J}(\tilde{y},\tilde{Z}_{.}) +C_{\alpha,K}|\ln \Delta t|^{-\frac{1}{2}}.\end{equation}
    Therefore, by Theorems \ref{thm-control} and \ref{thm-dis-error}, we have
    \begin{align}
        &\sup_{0\leq t\leq T}\E|X_{t}-\hat{X}_{t}|+\sup_{0\leq t\leq T}\E|Y_{t}-\hat{Y}_{t}|\nonumber\\
        \leq &\left[\sup_{0\leq t\leq T}\E|X_{t}-\tilde{X}_{t}|+\sup_{0\leq t\leq T}\E|Y_{t}-\tilde{Y}_{t}|\right]+\left[ \sup_{0\leq t\leq T}\E|\hat{X}_{t}-\tilde{X}_{t}|+\sup_{0\leq t\leq T}\E|\hat{Y}_{t}-\tilde{Y}_{t}|\right]\nonumber\\
        \leq & C\J(\tilde{y},\tilde{Z}_{.})+C_{\alpha,K}\left[ 1+\J(\tilde{y},\tilde{Z}_{.}) \right]|\ln \Delta t|^{-1}\nonumber\\
        \leq & C\hat{\J}(\tilde{y},\tilde{Z}_{.}) + C_{\alpha,K}\left[ |\ln \Delta t|^{-\frac{1}{2}}+\hat{\J}(\tilde{y},\tilde{Z}_{.}) |\ln \Delta t|^{-1}\right].
    \end{align}
    Furthermore, we have
    \begin{align}
        &\sup_{0\leq t\leq T}\E|X_{t}-\hat{X}_{t}|^{2}+\sup_{0\leq t\leq T}\E|Y_{t}-\hat{Y}_{t}|^{2}+\int_{0}^{T}\E|Z_{t}-\tilde{Z}_{t}|^{2}\dif t\nonumber\\
        \leq &C\left[\sup_{0\leq t\leq T}\E|X_{t}-\tilde{X}_{t}|^{2}+\sup_{0\leq t\leq T}\E|Y_{t}-\tilde{Y}_{t}|^{2}+\int_{0}^{T}\E|Z_{t}-\tilde{Z}_{t}|^{2}\dif t\right]\\
        &+C\left[ \sup_{0\leq t\leq T}\E|\hat{X}_{t}-\tilde{X}_{t}|^{2}+\sup_{0\leq t\leq T}\E|\hat{Y}_{t}-\tilde{Y}_{t}|^{2}\right]\nonumber\\
        \leq & C\left[\J(\tilde{y},\tilde{Z}_{.})+\J(\tilde{y},\tilde{Z}_{.})^{2}\right]+C_{\alpha,K}\left[ 1+\J(\tilde{y},\tilde{Z}_{.}) \right]|\ln \Delta t|^{-1}\nonumber\\
        \leq & C\left[\hat{\J}(\tilde{y},\tilde{Z}_{.})+\hat{\J}(\tilde{y},\tilde{Z}_{.})^{2}\right] +  C_{\alpha,K}\left[ |\ln \Delta t|^{-\frac{1}{2}}+\hat{\J}(\tilde{y},\tilde{Z}_{.}) |\ln \Delta t|^{-1}\right].
    \end{align}
\section{Numerical Results}
In this section, we present two numerical examples.
The first is a one-dimensional model for which we can obtain an analytical solution; the second is a multi-dimensional case that can only be solved numerically.
Recall that \(\tilde{Z}_{t_{k}}=\Phi(t_{k},\hat{X}_{t_{k}};\beta)\), where
\(\Phi\) is constructed with two hidden layers of dimensions \(n+10\).
We use rectifier function (ReLU) as the activation function and implement Adam optimizer.
Both examples in this section are computed based on 1000 sample paths.
The number of time partition is set to \(N=100\).
We run each example independently 10 times to get the average result.
All the parameters will be initialized using a uniform distribution.
Related source codes can be referred to \url{https://github.com/YifanJiang233/Deep_BSDE_solver}.

\begin{ex}
 The first example considers the bond price using the CIR model.
   Let \(X_{t}\) be the short rate following a CIR process
    \begin{equation*}
        \dif X_{t}= a(b-X_{t})\dif t+ \sigma \sqrt{X_{t}}\dif W_{t}.
    \end{equation*}
    Assume \(Y_{t}\) is the zero-coupon bond paying 1 at maturity \(T\).
    We have
    \begin{equation*}
        Y_{t}=\E \left[ \exp\left(-\int_{t}^{T}X_{s}\dif s\right)\Big| \mathcal{F}_{t} \right].
    \end{equation*}
    Actually, \(Y_{t}\) satisfies the following BSDE
    \begin{equation*}
        Y_{t}=1-\int_{t}^{T}X_{s}Y_{s}\dif s - \int_{t}^{T}Z_{s}\dif W_{s}.
    \end{equation*}
    On the other hand, \(Y_{t}\) has an explicit solution
    \begin{equation*}
        Y_{t}=\left( \frac{2\gamma\exp(\frac{(\gamma+a)(T-t)}{2})}{(\gamma-a)+(\gamma+a)\exp(\gamma(T-t))} \right)^{\frac{2ab}{\sigma^{2}}}\exp\left( \frac{2(1-\exp(\gamma(T-t)))X_{t}}{(\gamma-a)+(\gamma+a)\exp(\gamma(T-t))} \right),
    \end{equation*}
    where \(\gamma=\sqrt{a^{2}+2\sigma^{2}}\).
    It is straightforward to check that Assumptions 1-4 are satisfied.
\end{ex}
In this case, we have \(n=m=d=1\).
We set \(a=b=\sigma=T=X_{0}=1\), then from direct computation, the prescribed initial value is \(Y_{0}\doteq 0.39647\).
Initial parameter \(\alpha\) is selected by a uniform distribution on \([0,1]\).

\begin{figure}[H]
    \centering
    \includegraphics[width=\textwidth]{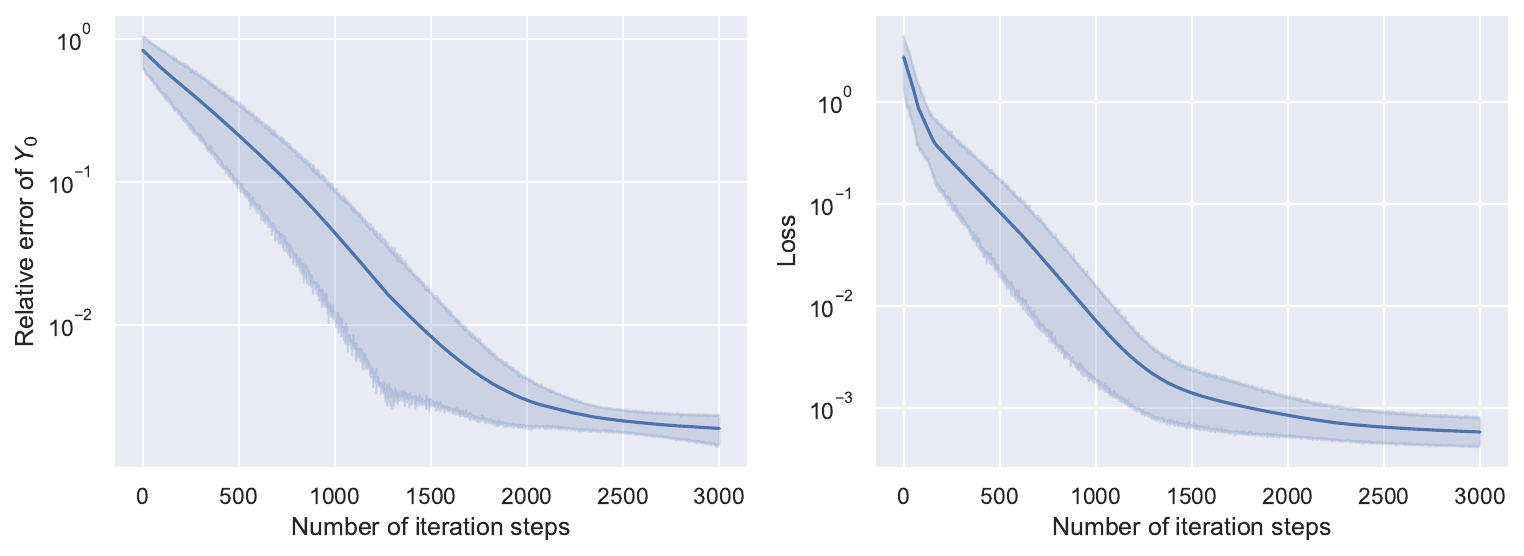}
    \caption{Relative error of the bond price (left) and the loss (right) against the number of the iteration steps.}
    \label{fig-cir}
\end{figure}

\begin{table}[H]
    \centering
    \caption{Numerical simulation of CIR bond}
    \begin{tabular}{|ccccc|}
        \hline Step & Mean of \(Y_{0}\) & Standard deviation of \(Y_{0}\) & Mean of loss & Standard deviation of loss \\
        \hline 500 & 0.4643            & 9.58E-2                         & 8.46E-2      & 1.27E-1                    \\
        \hline 1000 & 0.4136            & 2.55E-2                         & 7.13E-3      & 1.23E-2                    \\
        \hline 2000 & 0.3972            & 1.21E-3                         & 8.47E-4      & 6.23E-4                    \\
        \hline 3000 & 0.3972            & 3.69E-4                         & 5.80E-4      & 3.20E-4                    \\
        \hline
    \end{tabular}
\end{table}
We use algorithm 1 to obtain the numerical solution.
By means of Theorems \ref{thm-control} and \ref{thm-dis-error}, \(Y_{0}\) converges to the analytical solution.
In Figure \ref{fig-cir}, we illustrate the relative error of the bond price, and the loss is plotted against the number of iteration steps.
The shaded area represents the confidence interval after 10 independent runs.
Numerical simulations have well demonstrated the theoretical result.
\begin{ex}
    Our second example considers the bond price using the multi-dimensional CIR model.
    Let \(X_{t}\) be the short rate following a multi-dimensional CIR process, i.e., each component of \(X_{t}\) follows a CIR process,
    \begin{equation*}
        \dif X_{t}^{i}=a^{i}(b^{i}-X_{t}^{i})\dif t +\sigma^{i}\sqrt{X_{t}^{i}}\dif W_{t}.
    \end{equation*}
    Assume \(Y_{t}\) is the zero-coupon bond paying 1 at maturity \(T\).
    Under nonarbitrage condition, it is natural that
    \begin{equation*}
        Y_{t}=\E\left[ \exp\left( -\int_{t}^{T}\max_{1\leq i\leq n} X_{s}^{i}\dif s \right)\big| \mathcal{F}_{t} \right].
    \end{equation*}
    Hence, we have
    \begin{equation*}
        Y_{t}=1-\int_{t}^{T} \left( \max_{1\leq i\leq n}X^{i}_{s} \right)Y_{s}\dif s - \int_{t}^{T}Z_{s}\dif W_{s}.
    \end{equation*}
    Assume \(Y_{t}=u(t,X_{t})\).
    Following standard arguments, we have
    \begin{equation*}
        \left\{
        \begin{aligned}
             & u_{t}+\sum_{i=1}^{n}a^{i}(b^{i}-x^{i})\frac{\partial u}{\partial x^{i}}+\frac{1}{2}\sum_{1\leq i,\, j\leq n}\sigma^{i}\sigma^{j}\sqrt{x^{i}x^{j}}\frac{\partial^{2}u}{\partial x^{i}\partial x^{j}}- \left( \max_{1\leq i \leq n}x^{i} \right) u=0, \\
             & u(T,x)=1.
        \end{aligned}
        \right.
    \end{equation*}
\end{ex}
We set \(n=100\), \(m=d=T=1\), and \(X_{0}=(1,\cdots,1)\).
\(a^{i}\), \(b^{i}\), and \(\sigma^{i}\) are selected in \([0,1]\).
The initial value \(\alpha\) is selected by a uniform distribution on \([0,1]\).
Because solving the equation explicitly is difficult, we set the mean value of \(Y_{0}\) after 3000 iterations as the convergence limit and calculate the relative error of \(Y_{0}\) in 10 runs.

\begin{figure}[H]
    \centering
    \includegraphics[width=\textwidth]{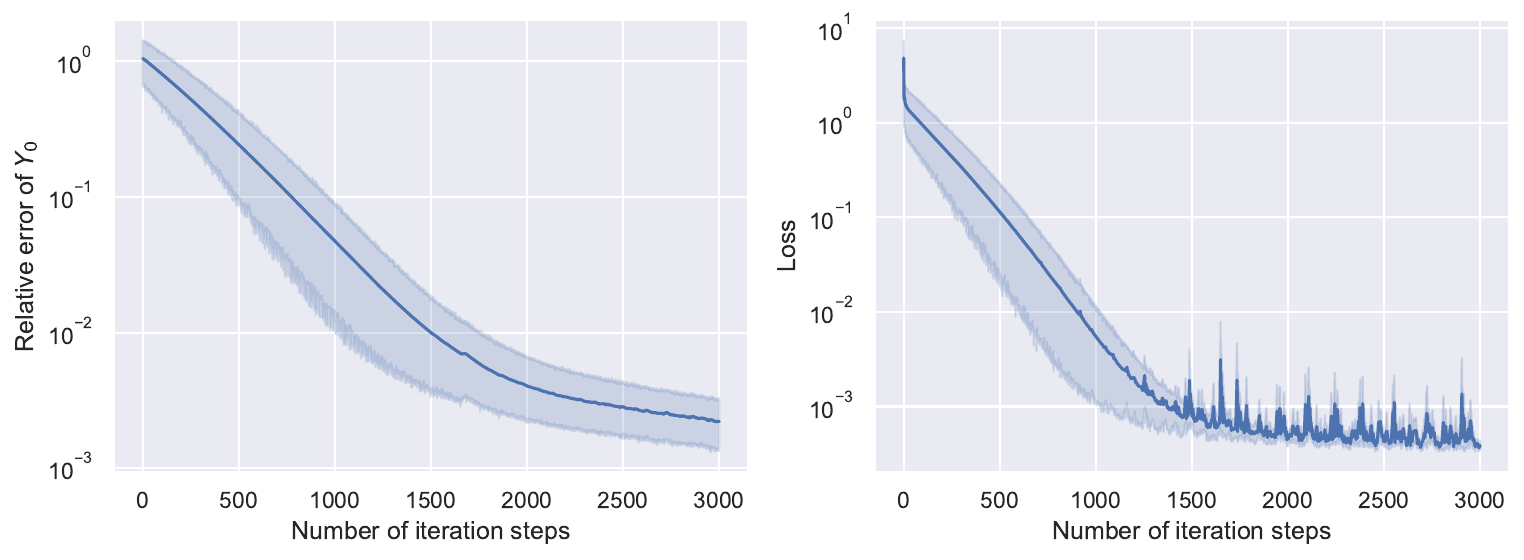}
    \caption{Relative error of the bond price under multi-dimensional CIR model (left) and the loss (right) against the number of the iteration steps.}
    \label{fig-multi-cir}
\end{figure}

In Figure \ref{fig-multi-cir}, we illustrate the relative error of the bond price under the multi-dimensional CIR model as well as a loss against the number of iteration steps.
The shaded area represents the confidence interval based on 10 independent runs. Numerical simulations show that \(Y_{0}\) also converges well after around \(2000\) iteration steps.

\begin{table}[H]
    \centering
    \caption{Numerical simulation of multi-dimensional CIR bond}
    \begin{tabular}{|ccccc|}
        \hline Step & Mean of \(Y_{0}\) & Standard deviation of \(Y_{0}\) & Mean of loss & Standard deviation of loss \\
        \hline 500 & 0.3773            & 8.77E-2                         & 1.15E-1      & 1.66E-1                    \\
        \hline 1000 & 0.3228            & 2.03E-2                         & 5.51E-3      & 8.33E-3                    \\
        \hline 2000 & 0.3100            & 1.63E-3                         & 4.50E-4      & 1.12E-4                    \\
        \hline 3000 & 0.3095            & 8.28E-4                         & 3.89E-4      & 7.20E-5                    \\
        \hline
    \end{tabular}
\end{table}
\section{Conclusions}
In this paper, we presented a posterior estimate to bound the error of a given numerical scheme for non-Lipschitz FBSDEs. 
We demonstrated that the proposed posterior estimate holds for continuous-time FBSDEs and extended it to the estimate for the discrete Euler-Maruyama scheme.
Some numerical examples in bond pricing are presented to demonstrate the Deep BSDE method's dependability.

We extend the results in previous studies \cite{han2020convergence,ji2020three} to a broader class of FBSDEs arising from the pricing problem in financial markets.
It is also worth mentioning that the posterior estimates in \cite{han2020convergence,ji2020three} require a sufficiently short time duration.
By adaptation of the four-step method \cite{ma1994solving}, we demonstrate the posterior estimate holds for arbitrary finite time duration if the decoupling field has bounded derivatives. 

Our next step is to extend the current result to McKean-Vlasov type FBSDEs where a careful analysis of the derivative of the decoupling field concerning distribution is required.
Another direction is to solve the Dirichlet problem of semilinear parabolic equations by the Deep BSDE method.
In this case, the first hitting time of the boundary is critical to the posterior estimate.

\bigskip

\section*{Acknowledgments}
This research has been supported by the EPSRC Centre for Doctoral Training in Mathematics of Random Systems: Analysis, Modelling, and Simulation (EP/S023925/1).

\bigskip
\bibliographystyle{amsplain}
\bibliography{mybib.bib}
\end{document}